\documentclass{article}
\usepackage{amsfonts}
\usepackage{amsmath}
\usepackage{amssymb}
\usepackage{graphicx}

\newcounter{guy}
\newtheorem{theorem}{Theorem}[section]
\newtheorem{lemma}{Lemma}[section]
\newtheorem{corollary}{Corollary}[section]

\newcommand{\ds}{\displaystyle}

\input{tcilatex}

\begin{document}

\title{Quadrangularity and Strong Quadrangularity in Tournaments}
\author{J. Richard Lundgren \\
\emph{\normalsize University of Colorado at Denver, }\\
\emph{\normalsize Denver, CO 80217, U.S.A }\\
\and K.B. Reid \\
\emph{\normalsize California State University San Marcos, }\\
\emph{\normalsize San Marcos, CA 92096, U.S.A.} \and Simone Severini \\
\emph{\normalsize University of York, }\\
\emph{\normalsize YO10 5DD, }\\
\emph{\normalsize Helsington, York, United Kingdom} \and Dustin J. Stewart 
\thanks{%
Corresponding author, e-mail: dstewart@math.cudenver.edu} \\
\emph{\normalsize University of Colorado at Denver, }\\
\emph{\normalsize Denver, CO 80217, U.S.A }}
\date{}
\maketitle

\begin{abstract}
The pattern of a matrix $M$ is a $(0,1)$-matrix which replaces all non-zero
entries of M with a 1. A directed graph is said to support $M$ if its
adjacency matrix is the pattern of $M$. If $M$ is an orthogonal matrix, then
a digraph which supports $M$ must satisfy a condition known as
quadrangularity. We look at quadrangularity in tournaments and determine for
which orders quadrangular tournaments exist. We also look at a more
restrictive necessary condition for a digraph to support an orthogonal
matrix, and give a construction for tournaments which meet this condition.
\end{abstract}

\section{Introduction}

A directed graph or digraph, $D,$ is a set of vertices $V(D)$ together with
a set of ordered pairs of the vertices, $A(D)$, called arcs. If $(u,v)$ is
an arc in a digraph, we say that $u$ beats $v$ or $u$ dominates $v$, and
typically write this as $u\rightarrow v$. If $v\in V(D)$ then we define the
outset of $v$ by, 
\begin{equation*}
O_{D}(v)=\{u\in V(D):(v,u)\in A(D)\}. 
\end{equation*}
That is, $O_{D}(v)$ is all vertices in $D$ which $v$ beats. Similarly, we
define the set of all vertices in $D$ which beat $v$ to be the inset of $v$,
written, 
\begin{equation*}
I_{D}(v)=\{u\in V(D):(u,v)\in A(D)\}. 
\end{equation*}
The closed outset and closed inset of a vertex $v$ are $O_{D}[v]=O_{D}(v)%
\cup\{v\}$ and $I_{D}[v]=I_{D}(v)\cup\{v\}$ respectively. The in-degree and
out-degree of a vertex $v$ are $d_{D}^{-}(v)=|I_{D}(v)|$ and $%
d_{D}^{+}(v)=|O_{D}(v)|$ respectively. When it is clear to which digraph $v$
belongs, we will drop the subscript. The minimum out-degree (in-degree) of $D
$ is the smallest out-degree (in-degree) of any vertex in $D$ and is
represented by $\delta^{+}(D)$ ($\delta^{-}(D)$). Similarly, the maximum
out-degree (in-degree) of $D$ is the largest out-degree (in-degree) of any
vertex in $D$ and is represented by $\Delta^{+}(D)$ ($\Delta^{-}(D)$).

\medskip\indent A tournament $T$ is a directed graph with the property that
for each pair of distinct vertices $u,v\in V(T)$ exactly one of $(u,v)$, $%
(v,u)$ is in $A(T)$. An $n$-tournament is a tournament on $n$ vertices. If $T
$ is a tournament and $W\subseteq V(T)$ we denote by $T[W]$ the
subtournament of $T$ induced on $W$. The dual of a tournament $T$, which we
denote by $T^{r}$, is the tournament on the same vertices as $T$ with $%
x\rightarrow y$ in $T^{r}$ if and only if $y\rightarrow x$ in $T$. If $%
X,Y\subseteq V(T)$ such that $x\rightarrow y$ for all $x\in X$ and $y\in Y$,
then we write $X\Rightarrow Y$. If $X=\{x\}$ or $Y=\{y\}$ we write $%
x\Rightarrow Y$ or $X\Rightarrow y$ respectively for $X\Rightarrow Y$. A
vertex $s\in V(T)$ such that $s\Rightarrow V(T)-s$ is called a transmitter.
Similarly a receiver is a vertex $t$ of $T$ such that $V(T)-t\Rightarrow t$.

\medskip\indent We say that a tournament is regular if every vertex has the
same out-degree. A tournament is called near regular if the largest
difference between the out-degrees of any two vertices is $1$. Let $S$ be a
subset of $\{1,2,\ldots,2k\}$ of order $k$ such that if $i,j\in S,$ $i+j\not
\equiv 0\pmod{2k+1}$. The tournament on $2k+1$ vertices labeled $%
0,1,\ldots,2k,$ with $i\rightarrow j$ if and only if $j-i\pmod{2k+1}\in S$
is called a rotational tournament with symbol $S$. If $p\equiv3\pmod{4}$ is
a prime and $S$ is the set of quadratic residues modulo $p$, then the
rotational tournament whose symbol is $S$ is called the quadratic residue
tournament of order $p$, denoted $QR_{p}$. We note that $|O(x)\cap
O(y)|=|I(x)\cap I(y)|=k$ for all distinct $x,y\in V(QR_{p})$ where $p=4k+3$.
For more on tournaments the reader is referred to \cite{Reid/Beineke}, \cite%
{Moon}, and \cite{ReidCRC}.

\medskip\indent Let $x=(x_{1},x_{2},\ldots,x_{n})$ and $y=(y_{1},y_{2},%
\ldots,y_{n})$ be $n$-vectors over some field (While the following
definitions hold over any field, we are interested only in those of
characteristic $0$). We use $\langle x,y\rangle$ to denote the usual
euclidean inner product of $x$ and $y$. We say that $x$ and $y$ are
combinatorially orthogonal if $|\{i:x_{i}y_{i}\neq0\}|\neq1$. Observe, this
is a necessary condition for $x$ and $y$ to be orthogonal, for if there were
a unique $i$ so that $x_{i}y_{i}\neq0$, then $\langle
x,y\rangle=x_{i}y_{i}\neq0$. We say a matrix $M$ is combinatorially
orthogonal if every two rows of $M$ are combinatorially orthogonal and every
two columns of $M$ are combinatorially orthogonal. In \cite{Klee}, Beasley,
Brualdi and Shader study matrices with the combinatorial orthogonality
property to obtain a lower bound on the number of non-zero entries in a
fully indecomposable orthogonal matrix.

\medskip\indent Let $M$ be an $n\times n$ matrix. The pattern of $M$ is the $%
(0,1)$-matrix whose $i,j$ entry is $1$ if and only if the $i,j$ entry of $M$
is non-zero. If $D$ is the directed graph whose adjacency matrix is the
pattern of $M$, we say that $D$ supports $M$ or that $D$ is the digraph of $M
$. We say a digraph $D$ is out-quadrangular if for all distinct $u,v\in V(D)$%
, $|O(u)\cap O(v)|\neq1$. Similarly, if for all distinct $u,v\in V(D)$, $%
|I(u)\cap I(v)|\neq1$, we say $D$ is in-quadrangular. If $D$ is both
out-quadrangular and in-quadrangular, then we say $D$ is quadrangular. It is
easy to see that if $D$ is the digraph of $M$, then $D$ is quadrangular if
and only if $M$ is combinatorially orthogonal. So, if $D$ is the digraph of
an orthogonal matrix, $D$ must be quadrangular. In \cite{Gibson}, Gibson and
Zhang study an equivalent version of quadrangularity in undirected graphs.
In \cite{Stew}, Lundgren, Severini and Stewart study quadrangularity in
tournaments. In the following section we expand on the results in \cite{Stew}%
, and in section $3$ we consider another necessary condition for a digraph
to support an orthogonal matrix.

\section{Known orders of quadrangular tournaments}

In this section we determine for exactly which $n$ there exists a
quadrangular tournament on $n$ vertices. We first need some results from 
\cite{Stew}.

\begin{theorem}
\cite{Stew} \label{odeg23} Let $T$ be an out-quadrangular tournament and
choose $v\in V(T)$. Let $W$ be the subtournament of $T$ induced on the
vertices of $O(v)$. Then $W$ contains no vertices of out-degree $1$.
\end{theorem}

\begin{theorem}
\cite{Stew} \label{ideg23} Let $T$ be an in-quadrangular tournament and
choose $v\in V(T)$. Let $W$ be the subtournament of $T$ induced on $I(v)$.
Then $W$ contains no vertices of in-degree $1$.
\end{theorem}

\begin{corollary}
\cite{Stew} \label{tourn8} If $T$ is an out-quadrangular tournament with $%
\delta^{+}(T)\geq2,$ then $\delta^{+}(T)\geq4.$
\end{corollary}

\begin{corollary}
\cite{Stew} \label{tourn8b} If $T$ is a quadrangular tournament with $%
\delta^{+}(T)\geq2$ and $\delta^{-}(T)\geq2$, then \newline
$\delta^{+}(T)\geq4$ and $\delta^{-}(T)\geq4.$
\end{corollary}

\medskip\indent Note that the only tournament on $4$ vertices with no vertex
of out-degree $1$ is a $3$-cycle together with a receiver. Similarly, the
only tournament on $4$ vertices with no vertex of in-degree $1$ is a $3$%
-cycle with a receiver. Thus, if a quadrangular tournament $T$ has a vertex $%
v$ of out-degree $4$, $T[O(v)]$ must be a $3$-cycle with a receiver, and if $%
u$ has in-degree $4$, $T[I(u)]$ must be a $3$-cycle with a transmitter.

\begin{theorem}
\label{order1} There does not exist a quadrangular near regular tournament
of order $10$.
\end{theorem}

\noindent Proof.~~ Suppose $T$ is such a tournament and pick a vertex $x$
with $d^{+}(x)=5$. So $d^{-}(x)=4$. Therefore $I(x)$ must induce a
subtournament comprised of a $3$-cycle, and a transmitter. Call this
transmitter $u$. If a vertex $y$ in $O(x)$ has $O(y)=I(x)$, then $|O(y)\cap
O(w)|=1$ for all $w\neq u$ in $I(x)$. This contradicts $T$ being
quadrangular, so $O(y)\neq I(x)$ for any $y\in O(x)$. Since every vertex in $%
O(x)$ beats at most $3$ vertices outside of $O(x)$, and since $T$ is near
regular we have that $\delta ^{+}(T[O(x)])\geq1$. Thus, by Theorem~\ref%
{odeg23}, we have $\delta ^{+}(T[O(x)])\geq2$. This means that $T[O(x)]$
must be the regular tournament on $5$ vertices.

\medskip\indent Consider the vertex $u$ which forms the transmitter in $%
T[I(x)]$. Since $u$ beats $I[x]-u$, and $T$ is near regular, $u$ can beat at
most one vertex in $O(x)$. If $u\rightarrow z$ for any $z\in O(x)$, then $%
|O(u)\cap O(x)|=|\{z\}|=1$ which contradicts $T$ being quadrangular. Thus, $%
z\rightarrow u$ for all $z\in O(x)$.

\medskip\indent Since $T$ is near regular, it has exactly $5$ vertices of
out-degree $5$, one of which is $x$. So, there can be at most four vertices
in $O(x)$ with out-degree $5$. Thus, there exists some vertex in $O(x)$ with
out-degree $4$, call it $v$. Since $x\rightarrow v$, $v$ beats $2$ vertices
in $O(x)$ and $v\rightarrow u$ there is exactly one vertex $r\in I(x)-u$
such that $v\rightarrow r$. Since $O(u)=I[x]-u$, we have $|O(v)\cap
O(u)|=|\{r\}|=1$. Therefore, $T$ is not quadrangular, and so such a
tournament does not exist. \hfill\framebox[.25cm]{~}

\medskip\indent Given a digraph $D$, and set $S\subseteq V(D)$, we say that $%
S$ is a dominating set in $D$ if each vertex of $D$ is in $S$ or dominated
by some vertex of $S$. The size of a smallest dominating set in $D$ is
called the domination number of $D$, and is denoted by $\gamma(D)$. In \cite%
{Stew} a relationship is shown to hold in certain tournaments between
quadrangularity and the domination number of a subtournament.

\begin{lemma}
\label{order2a} If $T$ is a tournament on $8$ vertices with $\gamma(T)\geq3$
and $\gamma(T^{r})\geq3$, then $T$ is near regular. Further, if $d^{-}(x)=3$%
, then $I(x)$ induces a $3$-cycle, and if $d^{+}(y)=3,$ then $O(y)$ induces
a $3$-cycle.
\end{lemma}

\noindent Proof.~~ Let $T$ be such a tournament. If $T$ has a vertex $a$
with $d^{-}(a)=0$ or $1$, then $I[a]$ would form a dominating set of size $1$
or $2$ respectively. If $T$ had a vertex $b$ with $I(b)=\{u,v\}$, where $%
u\rightarrow v$, then $\{u,b\}$ forms a dominating set of size $2$. So $%
d_{T}^{-}(x)\geq3$ for all $x\in V(T)$. Similarly, $d_{T^{r}}^{-}(x)\geq3$
for all $x\in V(T)$. Thus, 
\begin{equation*}
3\leq d_{T^{r}}^{-}(x)=d_{T}^{+}(x)=8-1-d_{T}^{-}(x)\leq7-3=4 
\end{equation*}
for all $x\in V(T)$. That is $3\leq d_{T}^{+}(x)\leq4$ for all $x\in V(T)$,
and $T$ is near regular. Now, pick $x\in V(T)$ with $d^{-}(x)=3$. If $I(x)$
induces a transitive triple with transmitter $u$, then $\{u,x\}$ would form
a dominating set in $T$. Thus, $I(x)$ must induce a $3$-cycle. By duality we
have that $O(y)$ induces a $3$-cycle for all $y$ with $d^{+}(y)=3$. \hfill%
\framebox[.25cm]{~}

\medskip\indent Up to isomorphism there are $4$ tournaments on $4$ vertices,
and exactly one of these is strongly connected. We refer to this tournament
as the strong $4$-tournament, and note that it is also the only tournament
on $4$ vertices without a vertex of out-degree $3$ or $0$.

\begin{lemma}
\label{order2b} Suppose $T$ is a tournament on $8$ vertices with $%
\gamma(T)\geq3$ and $\gamma(T^{r})\geq3$. Then if $x\in V(T)$ with $%
d^{+}(x)=4$, $O(x)$ induces the strong $4$-tournament.
\end{lemma}

\noindent Proof.~~ By Lemma~\ref{order2a}, $T$ is near regular so pick $x\in
V(T)$ with $d^{+}(x)=4$, and let $W$ be the subtournament induced on $O(x)$.
If there exists $u\in V(W)$ with $d_{W}^{+}(u)=0$, then since $%
d_{T}^{+}(u)\geq3$, $u\Rightarrow I(x)$ and $\{u,x\}$ forms a dominating set
in $T$. This contradicts $\gamma(T)\geq3$, so no such $u$ exists. Now assume
there exists a vertex $v\in V(W)$ with $d_{W}^{+}(v)=3$. If $d_{T}^{+}(v)=4$%
, then $v\rightarrow y$ for some $y\in I(x)$. So, $I(v)=I[x]-y$. However, $%
I(v)=I[x]-y$ forms a transitive triple, a contradiction to Lemma~\ref%
{order2a}. So $d_{T}^{+}(v)=3$. Now, since $\delta^{+}(W)>0$, the vertices
of $W-v$ all have out-degree $1$ in $W$. If some $z\in V(W)-v$ had $%
d_{T}^{+}(z)=4$, then $z\Rightarrow I(x)$ and $\{x,z\}$ would form a
dominating set of size $2$. Therefore, all $z\in V(W)$ have $d_{T}^{+}(z)=3$%
. Since $T$ is near regular, this implies that every vertex of $I[x]$ must
have out-degree $4$. Further, since $d_{T}^{+}(v)=3$, $O(v)\subseteq O(x)$
and so $I(x)\Rightarrow v$. So, each vertex of $I(x)$ dominates $x,v$ and
another vertex of $I(x)$. Thus, each vertex of $I(x)$ dominates a unique
vertex of $O(x)-v$. Further each vertex of $O(x)-v$ has out-degree $3$ in $T$
and so must be dominated by a unique vertex of $I(x)$. So label the vertices
of $I(x)$ as $y_{1},y_{2},y_{3}$ and the vertices of $O(x)-v$ as $%
w_{1},w_{2},w_{3}$ so that $y_{i}\rightarrow w_{i}$, and $w_{i}\rightarrow
y_{j}$ for $i\neq j$. Since $I(x)$ and $O(x)-v$ form $3$-cycles we may also
assume that $y_{1}\rightarrow y_{2}\rightarrow y_{3}$, $y_{3}\rightarrow
y_{1}$ and $w_{1}\rightarrow w_{2}\rightarrow w_{3}$ and $w_{3}\rightarrow
w_{1}$. So, $O(w_{1})=\{w_{2},y_{2},y_{3}\}$ which forms a transitive triple
a contradiction to Lemma~\ref{order2a}. Hence, no such $v$ exists and $%
1\leq\delta^{+}(W)\leq\Delta^{+}(W)\leq2$ and $W$ is the strong $4$%
-tournament. \hfill\framebox[.25cm]{~}

\begin{theorem}
\label{order2} Let $T$ be a tournament on $8$ vertices. Then $\gamma(T)\leq2$
or $\gamma(T^{r})\leq2$.
\end{theorem}

\noindent Proof.~~ Suppose to the contrary that $T$ is a tournament on $8$
vertices with $\gamma(T)\geq3$ and $\gamma(T^{r})\geq3$. By Lemma~\ref%
{order2a} we know that $T$ is near regular. Let $W$ be the subtournament of $%
T$ induced on the vertices of out-degree $4$. We can always choose $x$ in $W$
with $d_{W}^{-}(x)\geq2$. So pick $x\in V(T)$ with $d_{T}^{+}(x)=4$ so that
it dominates at most one vertex of out-degree $4$. By Lemma~\ref{order2b}, $%
O(x)$ induces the strong $4$-tournament. By our choice of $x$, at at least
one of the vertices with out-degree $2$ in $T[O(x)]$ has out-degree $3$ in $T
$. Call this vertex $x_{1}$. Label the vertices of $O(x_{1})\cap O(x)$ as $%
x_{2}$ and $x_{3}$ so that $x_{2}\rightarrow x_{3}$, and label the remaining
vertex of $O(x)$ as $x_{0}$. Note since $T[O(x)]$ is the strong $4$%
-tournament, we must have $x_{3}\rightarrow x_{0}$ and $x_{0}\rightarrow
x_{1}$. Since $d_{T}^{+}(x_{1})=3$, $x_{1}$ must dominate exactly one vertex
in $I(x)$, call it $y_{1}$. Recall $I(x)$ must induce a $3$-cycle by Lemma~%
\ref{order2a}, so we can label the remaining vertices of $I(x)$ as $y_{2}$
and $y_{3}$ so that $y_{1}\rightarrow y_{2}\rightarrow y_{3}$ and $%
y_{3}\rightarrow y_{1}$. Note since $O(x_{1})\cap I(x)=y_{1}$, $%
y_{2}\rightarrow x_{1}$ and $y_{3}\rightarrow x_{1}$. Also, by Lemma~\ref%
{order2a}, $O(x_{1})$ forms a $3$-cycle, so $x_{3}\rightarrow y_{1}$ and $%
y_{1}\rightarrow x_{2}$.

\medskip\indent Now, assume to the contrary that $y_{1}\rightarrow x_{0}$.
Then $O(y_{1})=\{x_{0},x_{2},x,y_{2}\}$. Now, since $O(x_{3})\cap
O(x)=\{x_{0}\}$, $d_{T}^{+}(x_{3})=3$ or else $x_{3}\Rightarrow I(x)$ and $%
\{x,x_{3}\}$ forms a dominating set of size $2$. So, $x_{3}$ dominates
exactly one of $y_{2}$ or $y_{3}$. If $x_{3}\rightarrow y_{2}$ then $%
y_{3}\rightarrow x_{3}$ and since $y_{3}\rightarrow x_{1}$, $\{y_{1},y_{3}\}$
forms a dominating set of size $2$. So, assume $x_{3}\rightarrow y_{3}$ and $%
y_{2}\rightarrow x_{3}$. Then $x,y_{3},x_{1},x_{3}\in O(y_{2})$ and $%
\{y_{2},y_{1}\}$ forms a dominating set of size $2$. Thus $x_{0}\rightarrow
y_{1}$.

\medskip\indent If $x_{3}\rightarrow y_{2}$, then $\{x_{3},x\}$ forms a
dominating set of size $2$, a contradiction. So, $y_{2}\rightarrow x_{3}$.
Now, if $x_{3}\rightarrow y_{3}$ then $O(x_{3})=\{y_{1},y_{3},x_{0}\}$.
However, $y_{3}\rightarrow y_{1}$ and $x_{0}\rightarrow y_{1}$ so $O(x_{3})$
forms a transitive triple, a contradiction to Lemma~\ref{order2a}. Thus $%
y_{3}\rightarrow x_{3}$. Since $d_{T}^{+}(y_{3})\leq4$ and $%
y_{1},x,x_{1},x_{3}\in O(y_{3})$, these are all the vertices in $O(y_{3})$.
So, $x_{0}\rightarrow y_{3}$.

\medskip\indent If $x_{0}\rightarrow y_{2}$ then $x_{0}\Rightarrow I(x)$ and 
$\{x,x_{0}\}$ form a dominating set of size $2$, so $y_{2}\rightarrow x_{0}$%
. So, $x_{0},y_{3},x\in O(y_{2})$ and $y_{1},x_{2},x_{3}\in O(x_{1})$, and
so $\{y_{2},x_{1}\}$ forms a dominating set of size $2$. Therefore, such a
tournament cannot exist. \hfill\framebox[.25cm]{~}

\begin{theorem}
\label{order3} No tournament $T$ on $9$ vertices with $\delta^{+}(T)\geq2$
is out-quadrangular.
\end{theorem}

\noindent Proof.~~ Suppose to the contrary $T$ is such a tournament. Since $T
$ is out-quadrangular, and $\delta^{+}(T)\geq2$, by Corollary~\ref{tourn8}, $%
\delta^{+}(T)\geq4$. Since the order of $T$ is $9$, this means $T$ must be
regular. Pick a vertex $x\in V(T)$. Then $O(x)$ must induce a subtournament
which is a $3$-cycle together with a receiver. Call the receiver of this
subtournament $y$. Since $T$ is regular, $d^{+}(y)=4$. Since $I(y)=O[x]-y$,
this means $O(y)=I(x)$. So, $O(y)=I(x)$ must induce a subtournament which is
a $3$-cycle together with a receiver vertex. Call this receiver $z$. Since $%
d^{+}(z)=4$, $y\rightarrow z$ and $I(x)-z$ dominate $z$, $O(z)=O[x]-y$. Now, 
$x\Rightarrow O(x)-y$ and $O(x)-y$ is a $3$-cycle so $T[O(z)]$ must contain
a vertex of out-degree $1$. Hence, by Theorem~\ref{odeg23}, $T$ is not
out-quadrangular. Thus, no such tournament exists.

\hfill\framebox[.25cm]{~}

\begin{corollary}
\label{order3a} No tournament $T$ on $9$ vertices with $\delta^{-}(T)\geq2$
is in-quadrangular.
\end{corollary}

\noindent Proof.~~ Let $T$ be a tournament on $9$ vertices with $\delta
^{-}(T)\geq2$. Then $T^{r}$ is not out-quadrangular by Theorem~\ref{order3}.
Thus $T$ is not in-quadrangular. \hfill\framebox[.25cm]{~}

\medskip\indent We now state a few more results from \cite{Stew}.

\begin{theorem}
\cite{Stew} \label{tourn7} Let $T$ be a tournament on $4$ or more vertices
with a vertex $x$ of out-degree $1$, say $x\rightarrow y$. Then, $T$ is
quadrangular if and only if 
\begin{list}{\arabic{guy}.}{\usecounter{guy}}
\item $O(y)=V(T)-\{x,y\},$
\item $\gamma(T-\{x,y\})>2$,
\item $\gamma((T-\{x,y\})^{r})>2$.
\end{list}
\end{theorem}

\begin{theorem}
\cite{Stew} \label{tourn2} Let $T$ be a tournament on $3$ or more vertices
with a transmitter $s$ and receiver $t$. Then $T$ is quadrangular if and
only if both $\gamma(T-\{s,t\})>2$ and $\gamma((T-\{s,t\})^{r})>2.$
\end{theorem}

\begin{theorem}
\cite{Stew}\label{tourn3} Let $T$ be a tournament with a transmitter $s$ and
no receiver. Then $T$ is quadrangular if and only if, $\gamma(T-s)>2$, $T-s$
is out-quadrangular, and $\delta^{+}(T-s)\geq2.$
\end{theorem}

\begin{corollary}
\cite{Stew} \label{tourn4} Let $T$ be a tournament with a receiver $t$ and
no transmitter. Then $T$ is quadrangular if and only if $\gamma((T-t)^{r})>2$%
, $T-t$ is in-quadrangular, and $\delta^{-}(T-t)\geq2.$
\end{corollary}

\begin{corollary}
\label{order4} No quadrangular tournament of order $10$ exists.
\end{corollary}

\noindent Proof.~~ By Corollaries~\ref{tourn8b} and \ref{tourn4}, and by
Theorems~\ref{tourn7}, \ref{tourn2} and~\ref{tourn3}, a quadrangular
tournament $T$ must satisfy one of the following. 
\begin{list}{\arabic{guy}.}{\usecounter{guy}}
\item $\delta^{+}(T)\geq 4$, and hence $T$ is near regular.
\item $T$ has a transmitter $s$ and receiver $t$ such that
$\gamma(T-\{s,t\})>2$ and $\gamma((T-\{s,t\})^{r})>2$.
\item $T$ contains an arc $(x,y)$ such that
$O(y)=I(x)=V(T)-\{x,y\}$ and $\gamma(T-\{x,y\})>2$ and
$\gamma((T-\{x,y\})^{r})>2$.
\item $T$ has a transmitter $s$ and $T-s$ is out-quadrangular
with $\delta^{+}(T-s)\geq 2$.
\item $T$ has a receiver $t$ and $T-t$ is in-quadrangular with
$\delta^{-}(T-t)\geq 2$.
\end{list}Note, Theorem~\ref{order1} implies that case 1 is impossible. If 2
or 3 were satisfied, then there would be a tournament on $8$ vertices such
that it and its dual have domination number at least $3$, which contradicts
Theorem~\ref{order2}. If 4 were satisfied, then $T-s$ would be of order $9$
and out-quadrangular, a contradiction to Theorem~\ref{order3}. Similarly, 5
contradicts Corollary~\ref{order3a}. Thus, no quadrangular tournament on $10$
vertices exists. \hfill\framebox[.25cm]{~}

\medskip\indent For the construction in Theorem~\ref{order6} we need the
following theorem from \cite{Stew}.

\begin{theorem}
\cite{Stew} \label{rotquad} Let $T$ be a rotational tournament on $n\geq5$
vertices, with symbol $S$. Then, $T$ is quadrangular if and only if for all
integers $m$ with $1\leq m\leq\frac{n-1}{2}$ there exist distinct subsets $%
\{i,j\},\{k,l\}\subseteq S$ such that $(i-j)\equiv(k-l)\equiv m\pmod{n}.$
\end{theorem}

\begin{theorem}
\label{order6} There exist quadrangular tournaments of order $11,12$ and $13$%
.
\end{theorem}

\noindent Proof.~~ Consider the quadratic residue tournament of order $11$, $%
QR_{11}$. For all $u,v\in V(QR_{11})$, recall that $|O(u)\cap
O(v)|=|I(u)\cap I(v)|=\frac{11-3}{4}=2$. Thus, $QR_{11}$ is quadrangular.
Further, this implies that for any two vertices $u,v\in V(QR_{11})$ there
exists a vertex which dominates both $u$ and $v$, so $\gamma(QR_{11})>2$.
Also, since $QR_{11}$ is regular, $\delta^{+}(QR_{11})=5\geq2$. Let $W$ be
the tournament formed by adding a transmitter to $QR_{11}$. Then by Theorem~%
\ref{tourn3}, $W$ is quadrangular.

\medskip\indent Now, let $T$ be the rotational tournament on $13$ vertices
with symbol $S=\{1,2,3,5,6,9\}$. The following table gives the subsets of $S$
which satisfy Theorem~\ref{rotquad}. Thus, $T$ is quadrangular. 
\begin{equation*}
\begin{array}{|c|c|}
\hline
\mbox{$m$} & \mbox{subsets} \\ \hline
1 & \{2,1\},\{3,2\} \\ \hline
2 & \{3,1\},\{5,3\} \\ \hline
3 & \{5,2\},\{6,3\} \\ \hline
4 & \{6,2\},\{9,5\} \\ \hline
5 & \{6,1\},\{1,9\} \\ \hline
6 & \{9,3\},\{2,9\} \\ \hline
\end{array}
\end{equation*}

\hfill\framebox[.25cm]{~}

\begin{theorem}
\label{order5} There exists a quadrangular tournament of order $14$.
\end{theorem}

\noindent Proof.~~ Construct $T$ of order $14$ in the following way. Start
with a set $V$ of $14$ distinct vertices. Partition $V$ into $7$ sets of
order $2$ labeled $V_{0},V_{1},V_{2},\ldots,V_{6}$. Each $V_{i}$ is to
induce the $2$-tournament, and $V_{i}\Rightarrow V_{j}$ if and only if $j-i%
\pmod{7}$ is one of $1,2,4$. We show that the resulting $14$-tournament, $T$%
, is quadrangular.

\medskip\indent Note that the condensation of $T$ on $V_{0},\ldots,V_{6}$ is
just the quadratic residue tournament on $7$ vertices, $QR_{7}$. Now, $QR_{7}
$ has the property that $|O(x)\cap O(y)|=1$ for all $x,y\in V(QR_{7})$.
Thus, if $u,v\in V(T)$ such that $u\in V_{i},$ $v\in V_{j}$ for $i\neq j,$ $%
|O(u)\cap O(v)|=2$. Further, since $QR_{7}$ is regular of degree $3$, if $%
u,v\in V(T)$ with $u,v\in V_{i}$ then $|O(u)\cap O(v)|=6$. Thus, $|O(u)\cap
O(v)|\neq1$ for all $u,v\in V(T)$, and so $T$ is out-quadrangular. Further,
since $QR_{7}$ is isomorphic to its dual, a similar argument shows that $T$
is in-quadrangular and hence quadrangular. \hfill\framebox[.25cm]{~}

\begin{theorem}
\label{order7} If $n\geq15,$ then there exists a quadrangular tournament on $%
n$ vertices.
\end{theorem}

\noindent Proof.~~ Pick $n\geq15.$ Let $a_{1},a_{2},a_{3},\ldots,a_{l}$ be a
sequence of at least $3$ integers such that $a_{i}\geq5$ for each $i,$ and $%
\displaystyle\sum_{i=1}^{l}a_{i}=n.$ Pick $l$ regular or near regular
tournaments $T_{1},T_{2},\ldots,T_{l}$ such that $|V(T_{i})|=a_{i}$ for each 
$i.$ Let $T^{\prime}$ be a tournament with $V(T^{\prime})=\{1,2,3,\ldots l\}$
such that $T^{\prime}$ has no transmitter or receiver. Construct the
tournament $T$ on $n$ vertices as follows. Start with a a set $V$ of $n$
vertices, and partition $V$ into sets $S_{1},S_{2},\ldots,S_{l}$ of size $%
a_{1},a_{2},\ldots,a_{l}$ respectively. Place arcs in each $S_{i}$ to form $%
T_{i}.$ Now, add arcs such that $S_{i}\Rightarrow S_{j}$ if and only if $%
i\rightarrow j$ in $T^{\prime}.$ We claim that the resulting tournament, $T$%
, is quadrangular.

\medskip\indent Pick $u,v\in V(T)$. We consider two possibilities. First,
suppose that $u,v\in S_{i}$ for some $i.$ By choice of $T^{\prime},$ $%
i\rightarrow j$ for some $j.$ Thus 
\begin{equation*}
|O(u)\cap O(v)|\geq|S_{j}|=a_{j}\geq5>1. 
\end{equation*}
Now, suppose that $u\in S_{i}$ and $v\in S_{j}$ for $i\neq j.$ Since $%
T^{\prime}$ is a tournament either $i\rightarrow j$ or $j\rightarrow i.$
Without loss of generality, assume that $i\rightarrow j.$ Then 
\begin{equation*}
|O(u)\cap O(v)|\geq|O(v)\cap S_{j}|\geq\frac{k-1}{2}\geq2>1. 
\end{equation*}
This shows that $T$ is out-quadrangular. The proof that $T$ is
in-quadrangular is similar. Thus, $T$ is a quadrangular tournament of order $%
n\geq15.$ \hfill\framebox[.25cm]{~}

\medskip\indent Observe that if $T^{\prime}$ in the construction is strong,
then $T$ is strong. Further, if $a_{i}=k$ for all $i$ and $T^{\prime}$ is
regular, then $T$ is regular or near regular depending on if $k$ is odd or
even. We now characterize those $n$ for which there exist a quadrangular
tournament of order $n$.

\begin{theorem}
\label{order8} There exists a quadrangular tournament of order $n$ if and
only if $n=1,2,3,9$ or $n\geq11$.
\end{theorem}

\noindent Proof.~~ Note that the single vertex, the single arc, and the $3$%
-cycle are all quadrangular. Now, recall that the smallest tournament with
domination number $3$ is $QR_{7}$ (For a proof of this see \cite{Reid/McRae}%
). Further, $QR_{7}$ is isomorphic to its dual, so $\gamma(QR_{7}^{r})=3$.
This fact together with Theorems~\ref{tourn7} and \ref{tourn2} tell us that
the smallest quadrangular tournament, $T$, on $n\geq4$ vertices with $\delta
^{+}(T)=\delta^{-}(T)=0$ or $\delta^{+}(T)=1$ or $\delta^{-}(T)=1$ has order 
$9$.

\medskip\indent Theorem~\ref{tourn3} and Corollary~\ref{tourn4} together
with the fact that $QR_{7}$ is the smallest tournament with domination
number $3$ imply that a quadrangular tournament with just a transmitter or
receiver must have at least $8$ vertices. However, $QR_{7}$ is the only
tournament on $7$ vertices with domination number $3$ and a quick check
shows that $QR_{7}$ is neither out-quadrangular nor in-quadrangular. So, $%
QR_{7}$ together with a transmitter or receiver is not quadrangular, and
hence any quadrangular tournament with just a transmitter or receiver must
have order $9$ or higher.

\medskip\indent Corollary~\ref{tourn8b} states that if $\delta^{+}(T)\geq2$
and $\delta^{-}(T)\geq2$, then $\delta^{+}(T)\geq4$ and $\delta^{-}(T)\geq4$%
. The smallest tournament which meets these requirements is a regular
tournament on $9$ vertices. Thus, there are no quadrangular tournaments of
order $4,5,6,7$ or $8$. The result now follows from Corollary~\ref{order4}
and Theorems~\ref{order6}, \ref{order5} and \ref{order7}. \hfill %
\framebox[.25cm]{~}

\medskip\indent It turns out that quadrangularity is a common (asymptotic)
property in tournaments as the following probabilistic result shows.

\begin{theorem}
\label{prob} Almost all tournaments are quadrangular.
\end{theorem}

\noindent Proof.~~ Let $P(n)$ denote the probability that a random
tournament on $n$ vertices contains a pair of distinct vertices $x$ and $y$
so that $|O(x)\cap O(y)|=1$. We now give an over-count for the number of
labeled tournaments on $n$ vertices which contain such a pair, and show $%
P(n)\rightarrow0$ as $n\rightarrow\infty$.

\medskip\indent There are $\binom{n}{2}$ ways to pick distinct vertices $x$
and $y$, and the arc between them can be oriented so that $x\rightarrow y$
or $y\rightarrow x$. There are $n-2$ vertices which can play the role of $z$
where $\{z\}=O(x)\cap O(y)$. For each $w\not \in \{x,y,z\}$ there are $3$
ways to orient the arcs from $x$ and $y$ to $w$, namely $w\Rightarrow x,y$, $%
w\rightarrow x$ and $y\rightarrow w$, or $w\rightarrow y$ and $x\rightarrow w
$. Also, there are $n-3$ such $w$. The arcs between all other vertices are
arbitrary. So there are $2^{\binom{n-2}{2}}$ ways to finish the tournament.
When orienting the remaining arcs we may double count some of these
tournaments, so all together there are at most 
\begin{equation*}
2\binom{n}{2}(n-2)3^{n-3}2^{\binom{n-2}{2}}
\end{equation*}
tournaments containing such a pair of vertices. Now, there are $2^{\binom {n%
}{2}}$ total labeled tournaments so, 
\begin{align*}
0\leq P(n) & \leq\frac{2\binom{n}{2}(n-2)3^{(n-3)}2^{\binom{n-2}{2}}}{2^{%
\binom{n}{2}}} \\
& =\frac{n(n-1)(n-2)3^{(n-3)}2^{\binom{n-2}{2}}}{2^{\binom{n-2}{2}+n-2+n-1}}
\\
& =\frac{n(n-1)(n-2)3^{n-3}}{2^{2n-3}} \\
& =\frac{n(n-1)(n-2)3^{n-3}}{2^{2(n-3)}2^{3}} \\
& =\frac{n(n-1)(n-2)}{8}\left( \frac{3}{4}\right) ^{n-3} \\
& =\frac{\frac{1}{8}n(n-1)(n-2)}{(\frac{4}{3})^{n-3}}.
\end{align*}
Since this value tends to $0$ as $n$ tends to $\infty$, it must be that $%
P(n)\rightarrow0$ as $n\rightarrow\infty$.

\medskip\indent From duality we have that the probability that vertices $x$
and $y$ exists such that $|I(x)\cap I(y)|=1$ also tends to $0$ as $n$ tends
to $\infty$. Thus, the probability that a tournament is not quadrangular
tends to $0$ as $n$ tends to $\infty$. That is, almost all tournaments are
quadrangular. \hfill\framebox[.25cm]{~}

\section{Strong Quadrangularity}

In this section we define a stronger necessary condition for a digraph to
support an orthogonal matrix, and give a construction for a class of
tournaments which satisfy this condition. Let $D$ be a digraph. Let $%
S\subseteq V(D)$ such that for all $u\in S$, there exists $v\in S$ such that 
$O(u)\cap O(v)\neq\emptyset$, and let $S^{\prime}\subseteq V(D)$ such that
for all $u\in S^{\prime}$, there exists $v\in S^{\prime}$ such that $%
I(u)\cap I(v)\neq\emptyset$. We say that $D$ is strongly quadrangular if for
all such sets $S$ and $S^{\prime}$, 
\begin{list}{(\roman{guy})}{\usecounter{guy}}
	\item $\ds \left|\bigcup_{u,v\in S}(O(u)\cap
	O(v))\right|\geq|S|,$
	\item $\ds \left|\bigcup_{u,v\in S'}(I(u)\cap
	I(v))\right|\geq|S'|.$
\end{list}In \cite{Severini}, Severini showed that strong quadrangularity is
a necessary condition for a digraph to support an orthogonal matrix. To see
that this is in fact a more restrictive condition consider the following
tournament. Let $T$ be a tournament with $V(T)=\{0,1,2,3,4,5,6,x,y\}$ so
that $\{0,1,2,3,4,5,6\}$ induce the tournament $QR_{7}$, $x\rightarrow y$
and $O(y)=I(x)=V(T)-\{x,y\}$. In the previous section we saw that $T$ is
quadrangular. Now consider the set of vertices $S=\{0,1,5\}$. Since each of $%
0,1,5$ beat $x$, we have that for all $u\in S$, there exits $v\in S$ so that 
$O(u)\cap O(v)\neq\emptyset$. Also, 
\begin{align*}
\left\vert \bigcup_{u,v\in S}(O(u)\cap O(v))\right\vert & =\left\vert
(O(0)\cap O(1))\cup(O(0)\cap O(5))\cup(O(1)\cap O(5))\right\vert \\
& =\left\vert \{2,x\}\cup\{2,x\}\cup\{2,x\}\right\vert \\
& =2 \\
& <|S|.
\end{align*}
So $T$ is not strongly quadrangular. We now construct a class of strongly
quadrangular tournaments, but first observe the following lemma.

\begin{lemma}
\label{comp} Let $T$ be a tournament on $n\geq4$ vertices. Then there must
exist distinct $a,b\in V(T)$ such that $O(a)\cap O(b)\neq\emptyset$.
\end{lemma}

\noindent Proof.~~ Pick a vertex $a$ of maximum out-degree in $T$. As, $%
n\geq4$, $d^{+}(a)\geq2$. Pick a vertex $b$ of maximum out-degree in the
subtournament $W$ induced on $O(a)$. As $d^{+}(a)\geq2$, $d_{W}^{+}(b)\geq1$%
. Thus, $|O(a)\cap O(b)|=d_{W}^{+}(b)\geq1$. \hfill\framebox[.25cm]{~}

\begin{theorem}
Pick $l\geq1$. Let $T^{\prime}$ be a strong tournament on the vertices $%
\{1,2,\ldots,l\},$ and let $T_{1},T_{2},\ldots,T_{l}$ be regular or
near-regular tournaments of order $k\geq5.$ Construct a tournament $T$ on $kl
$ vertices as follows. Let $V$ be a set of $kl$ vertices. Partition the
vertices of $V$ into $l$ subsets $V_{1},\ldots,V_{l}$ of size $k$ and place
arcs to form copies of $T_{1},T_{2},\ldots,T_{l}$ on $V_{1},\ldots,V_{l}$
respectively. Finally, add arcs so that $V_{i}\Rightarrow V_{j}$ if and only
if $i\rightarrow j$ in $T^{\prime}.$ Then the resulting tournament, $T$, is
a strongly quadrangular tournament.
\end{theorem}

\noindent Proof.~~ Pick $S\subseteq V(T).$ Define the set 
\begin{equation*}
A=\{V_{i}:\exists u\neq v\in S\ni u,v\in V_{i}\}, 
\end{equation*}
and define the set 
\begin{equation*}
B=\{V_{i}:\exists!u\in S\ni u\in V_{i}\}. 
\end{equation*}
Let $\alpha=|A|,$ and $\beta=|B|.$ Then, since each $V_{i}$ has $k$
vertices, $k\alpha+\beta\geq|S|.$ Consider the subtournaments of $T^{\prime}$
induced on the vertices corresponding to $A$ and $B.$ These are tournaments
and so must contain a Hamiltonian path. So, label the elements of $A$ and $B$
so that $A_{1}\Rightarrow A_{2}\Rightarrow\cdots\Rightarrow A_{\alpha}$ and $%
B_{1}\Rightarrow B_{2}\Rightarrow\cdots\Rightarrow B_{\beta}.$ By definition
of $A,$ each $A_{i}$ contains at least two vertices of $S$, and so if $%
x,y\in S$ and $x,y\in A_{i},$ $i\leq\alpha-1,$ then $A_{i+1}\subseteq
O(x)\cap O(y).$ Thus, 
\begin{equation*}
\left\vert \bigcup_{u,v\in S}O(u)\cap O(v)\right\vert \geq k(\alpha-1). 
\end{equation*}
We now consider three cases depending on $\beta.$

\medskip\indent First assume that $\beta\geq2.$ Consider the vertices of $S$
in $B$ we see that if $x,y\in S$ and $x\in B_{i}$ and $y\in B_{i+1}$ then $%
O(y)\cap B_{i+1}\subseteq O(x)\cap O(y).$ Thus, $|O(x)\cap O(y)|\geq \frac{%
k-1}{2},$ and so 
\begin{equation*}
\left\vert \bigcup_{u,v\in S}O(u)\cap O(v)\right\vert \geq k(\alpha -1)+%
\frac{k-1}{2}(\beta-1)\geq k(\alpha-1)+2\beta-2\geq k(\alpha-1)+\beta. 
\end{equation*}
Now, since $T^{\prime}$ is a tournament, either $A_{1}\Rightarrow B_{1}$ or $%
B_{1}\Rightarrow A_{1}.$ If $A_{1}\Rightarrow B_{1},$ then for vertices $%
x,y\in A_{1}$ we know $B_{1}\subseteq O(x)\cap O(y).$ Since no vertex of $%
B_{1}$ had been previously counted, we have that 
\begin{equation*}
\left\vert \bigcup_{u,v\in S}O(u)\cap O(v)\right\vert \geq k(\alpha
-1)+\beta+k=k\alpha+\beta\geq|S|. 
\end{equation*}
So, assume that $B_{1}\Rightarrow A_{1}.$ Then for the single vertex of $S$
in $B_{1},$ $u,$ and a vertex $v$ of $S$ in $A_{1}$ $O(v)\subseteq O(u)\cap
O(v).$ This adds $\frac{k-1}{2}$ vertices which were not previously counted.
Also, since $T^{\prime}$ is strong, some $A_{i}\Rightarrow V_{j}$ for some $%
V_{j}\not \in A.$ We counted at most $\frac{k-1}{2}$ vertices in $V_{j}$
before, and since $A_{i}$ contains at least two vertices $x,y$ from $S$
these vertices add at least $\frac{k+1}{2}$ vertices which were not
previously counted, so 
\begin{equation*}
\left\vert \bigcup_{u,v\in S}O(u)\cap O(v)\right\vert \geq k(\alpha
-1)+\beta+\frac{k-1}{2}+\frac{k+1}{2}=k\alpha+\beta\geq|S|. 
\end{equation*}

\medskip\indent Now assume that $\beta=1.$ Since $T^{\prime}$ is strong we
know that $A_{i}\Rightarrow V_{j}$ for some $V_{j}\not \in A.$ So, 
\begin{equation*}
\left\vert \bigcup_{u,v\in S}O(u)\cap O(v)\right\vert \geq k\alpha. 
\end{equation*}
Now, if $|S|\leq k\alpha,$ then we are done, so assume that $|S|=k\alpha+1.$
So, for every $A_{i}\in A,$ $A_{i}\subseteq S.$ So by Lemma~\ref{comp} we
can find two vertices of $S$ in $A_{1}$ which compete over a vertex of $%
A_{i},$ adding one more vertex to our count, and 
\begin{equation*}
\left\vert \bigcup_{u,v\in S}O(u)\cap O(v)\right\vert \geq k\alpha+1\geq|S|. 
\end{equation*}

\medskip\indent For the last case, assume that $\beta=0.$ Then since $%
T^{\prime}$ is strong we once again have that some $A_{i}\Rightarrow V_{j}$
for some $V_{j}\not \in A.$ Thus, 
\begin{equation*}
\left\vert \bigcup_{u,v\in S}O(u)\cap O(v)\right\vert \geq k\alpha\geq|S|. 
\end{equation*}

\medskip\indent Note that the dual of $T^{\prime}$ will again be strong, and
the dual of each $T_{i}$ will again be regular. Thus, by appealing to
duality in $T$ we have that for all $S\subseteq V(T)$, 
\begin{equation*}
\left\vert \bigcup_{u,v\in S}I(u)\cap I(v)\right\vert \geq|S|, 
\end{equation*}
and so $T$ is a strongly quadrangular tournament.

\hfill\framebox[.25cm]{~}

\medskip\indent Recall that strong quadrangularity is a necessary condition
for a digraph to support an orthogonal matrix. To emphasize this, consider
the strongly quadrangular tournament, $T$, which the construction in the
previous theorem gives on $15$ vertices. For this tournament, $T_{1},T_{2}$
and $T_{3}$ are all regular of order $5$, and $T^{\prime}$ is the $3$-cycle.
Note that up to isomorphism, there is only one regular tournament on $5$
vertices, so without loss of generality, assume that $T_{1},T_{2}$ and $T_{3}
$ are the rotational tournament with symbol $\{1,2\}$. We now show that $T$
cannot be the digraph of an orthogonal matrix.

\medskip\indent Let $J_{5}$ denote the $5\times5$ matrix of all $1$s, $O_{5}$
the $5\times5$ matrix of all $0$s and set 
\begin{equation*}
RT_{5}=\left( 
\begin{array}{ccccc}
0 & 1 & 1 & 0 & 0 \\ 
0 & 0 & 1 & 1 & 0 \\ 
0 & 0 & 0 & 1 & 1 \\ 
1 & 0 & 0 & 0 & 1 \\ 
1 & 1 & 0 & 0 & 0%
\end{array}
\right) . 
\end{equation*}
Then the adjacency matrix $M$ of $T$ is 
\begin{equation*}
M=\left( 
\begin{array}{ccc}
RT_{5} & J_{5} & O_{5} \\ 
O_{5} & RT_{5} & J_{5} \\ 
J_{5} & O_{5} & RT_{5}%
\end{array}
\right) . 
\end{equation*}
Now, suppose to the contrary that there exists an orthogonal matrix $U$
whose pattern is $M$. Let $R_{i}$ and $C_{i}$ denote the $i^{th}$ rows and
columns of $U$ respectively for each $i=1,\ldots,15$, and let $U_{i,j}$
denote the $i,j$ entry of $U$. Observe from the pattern of $U$ that the only
entries of $U$ which contribute to $\langle C_{i},C_{j}\rangle$ for $%
i=1,\ldots,5$, $j=6,\ldots,10$ are in the first five rows. So, $\langle
C_{1},C_{j}\rangle=U_{4,1}U_{4,j}+U_{5,1}U_{5,j}$ for $j=6,\ldots,10$. Thus,
since $0=\langle C_{1},C_{j}\rangle$ for each $j\neq1$, 
\begin{equation*}
U_{4,1}=\frac{-U_{5,1}U_{5,6}}{U_{4,6}}=\frac{-U_{5,1}U_{5,7}}{U_{4,7}}=%
\frac{-U_{5,1}U_{5,8}}{U_{4,8}}=\frac{-U_{5,1}U_{5,9}}{U_{4,9}}=\frac{%
-U_{5,1}U_{5,10}}{U_{4,10}}. 
\end{equation*}
Since $U_{5,1}\neq0$ this gives, 
\begin{equation*}
-\frac{U_{4,1}}{U_{5,1}}=\frac{U_{5,6}}{U_{4,6}}=\frac{U_{5,7}}{U_{4,7}}=%
\frac{U_{5,8}}{U_{4,8}}=\frac{U_{5,9}}{U_{4,9}}=\frac{U_{5,10}}{U_{4,10}}. 
\end{equation*}
So, the vectors $(U_{4,6},\ldots,U_{4,10})$ and $(U_{5,6},\ldots,U_{5,10})$
are scalar multiples of each other. Now, note that for $j=6,\ldots,10$, we
have $0=\langle C_{2},C_{j}\rangle=U_{1,2}U_{1,j}+U_{5,2}U_{5,j}$. So, by
applying the same argument, we see that $%
(U_{5,6},U_{5,7},U_{5,8},U_{5,9},U_{5,10})$ is a scalar multiple of $%
(U_{1,6},U_{1,7},U_{1,8},U_{1,9},U_{1,10})$. So, $%
(U_{4,6},U_{4,7},U_{4,8},U_{4,9},U_{4,10})$ is a scalar multiple of $%
(U_{1,6},U_{1,7},U_{1,8},U_{1,9},U_{1,10})$. Now, from the pattern of $U$ we
see that only the $6^{th}$ through $10^{th}$ columns of $U$ contribute to $%
\langle R_{1},R_{4}\rangle$. So, since linearly dependent vectors cannot be
orthogonal, 
\begin{equation*}
\langle
R_{1},R_{4}\rangle=%
\langle(U_{1,6},U_{1,7},U_{1,8},U_{1,9},U_{1,10}),(U_{4,6},U_{4,7},U_{4,8},U_{4,9},U_{4,10})\rangle\neq0. 
\end{equation*}
This contradicts our assumption that $U$ is orthogonal. So, $T$ is not the
digraph of an orthogonal matrix.

\section{Conclusions}

The problem of determining whether or not there exist tournaments (other
than the $3$-cycle) which support orthogonal matrices has proved to be quite
difficult. As we have seen in sections 2 and 3, for large values of $n$ we
can almost always construct examples of tournaments which meet our necessary
conditions. Knowing that almost all tournaments are quadrangular and having
a construction for an infinite class of strongly quadrangular tournaments,
one may believe that there will exist a tournament which supports an
orthogonal matrix. However, attempting to find an orthogonal matrix whose
digraph is a given tournament has proved to be a difficult task. In general,
aside from the $3$-cycle, the existence of a tournament which supports an
orthogonal matrix is still an open problem. We conclude this section with a
result that may lead one to believe non-existence is the answer to this
problem.

\begin{theorem}
Other than the $3$-cycle, there does not exist a tournament on $10$ or fewer
vertices which is the digraph of an orthogonal matrix.
\end{theorem}

\noindent Proof.~~ By Theorem~\ref{order8}, there exists a quadrangular $n$%
-tournament for $n\leq10$ if and only if $n$ is $1$, $2$, $3$ or $9$. Note,
in the case $n=1$ and $n=2$, the only tournaments are the single vertex and
single arc, both of whose adjacency matrices have a column of zeros. Since
orthogonal matrices have full rank, these cannot support an orthogonal
matrix. When $n=3$, the $3$-cycle is the only quadrangular tournament. The
adjacency matrix for this tournament is a permutation matrix and hence
orthogonal. Now consider $n=9$. By Theorem~\ref{order3}, if $T$ is
quadrangular, $\delta ^{+}(T)\leq1$. If $\delta^{+}(T)=0$, then $T$'s
adjacency matrix will have a row of zeros, and $T$ cannot be the digraph of
an orthogonal matrix. So we must have $\delta^{+}(T)=1$. So by Theorem~\ref%
{tourn7}, $T$ has an arc $(x,y)$ with $O(y)=I(x)=V(T)-\{x,y\}$ and $%
\gamma(T-\{x,y\})>2$. The only $7$-tournament with domination number greater
than $2$ is $QR_{7}$, thus $T-\{x,y\}=QR_{7}$. However, in section 3 we
observed that this tournament is not strongly quadrangular. Thus, other than
the $3$-cycle, no tournament on $10$ or fewer vertices can be the digraph of
an orthogonal matrix. \hfill\framebox[.25cm]{~}

\end{document}